\documentclass{amsart}
\usepackage{amsfonts}

\usepackage{amsmath}
\usepackage{graphicx}
\usepackage{amscd}

\newtheorem{theorem}{Theorem}
\theoremstyle{plain}

\newtheorem{remark}{Remark}

\numberwithin{equation}{section}

\input{tcilatex}

\begin{document}
\title[On Barnes' type multiple Changhee
$q$-Zeta functions]{ Barnes' type multiple  Changhee $q$-Zeta
functions}
\author{TAEKYUN KIM}
\address{TAEKYUN KIM\\
Institute of Science Education, Kongju National University Kongju 314-701,
S. Korea}
\email{tkim64@hanmail.net}
\author{YILMAZ\ SIMSEK}
\curraddr{YILMAZ\ SIMSEK\\
Mersin University, Faculty of Science, Department of Mathematics 33343
Mersin, Turkey}
\email{ysimsek@mersin.edu.tr}

\begin{abstract}
In this paper, we give new generating functions which produce
Barnes' type multiple generalized Changhee $q$-Bernoulli
polynomials and polynomials. These functions are very important to
construct multiple zeta functions. By using Mellin transform's
formula and Cauchy Theorem , we prove the analytic continuation of
Barnes' type multiple Changhee $q$-zeta function. Finally we give
some relations between Barnes' type multiple Changhee $q$-zeta
function and Barnes' type multiple\ generalized Changhee
$q$-Bernoulli numbers.
\end{abstract}

\maketitle

\section{Introduction, Definition and Notations}

Multiple zeta functions have studied by many mathematicians. These functions
are of interest and importance in many areas. These functions and numbers
are in used not only in Complex Analysis and Mathematical Physics, but also
in used in $p$-adic Analysis and other areas. In particular, multiple zeta
functions occur within the context of knot theory, Quantum Field Theory,
Applied Analysis and Number TheoryIn particular, multiple zeta functions
occur within the context of knot theory and quantum field theory. The Barnes
multiple zeta functions and gamma functions were also encountered by
Shintani within the context of analytic number theory. They showed up in the
form of the factor program for integrable field theories and in the studies
of $XXZ$ model correlation function. The computation of Feynman diagrams has
confronted physicists with classes of integrals that are usually hard to be
evaluated, both analytically and numerically. The newer techniques applied
in the more popular computer algebra packages do not offer much relief.
Therefore, it seems reasonable to occasionally study some alternative
methods to come to a result. In the case of the computation of structure
functions in deep inelastic scattering, the Mellin moments of these
functions are often of interest. Each individual moment can be computed
directly in a much simples way than that needed to compute the entire
structure function and take its moments afterwards. Meanwhile, the special
values of multiple zeta functions at positive integers have come to the
foreground in the recent years, both in connection with theorical physics
(Feynman diagrams) and the theory of mixed Tate motives. Historically, Euler
already investigated the double zeta values in the $XVIII$th century ( see\
for detail \cite{Andrews}, \cite{Askey}, \cite{Barnes}, \cite{Cherednik}, %
\cite{C. A. Nelson and M. G. Gartley}, \cite{B. E. Sagan}, \cite{T. H.
Koornwinder}, \cite{A. M. Robert}, \cite{E. T. Wittaker and G. N. Watson}, %
\cite{Washington},\ \cite{Katriel}, \cite{Khrennikov}, \cite{M. Nishizawa},
\cite{N. Kurokawa}, \cite{W. H. Schikhof}, \cite{Kim1}, \cite{Kim2}, \cite%
{Kim4}, \cite{Kim5}, \cite{Kim7}, \cite{Kim11}, \cite{Kim14}, \cite%
{Kim-Simsek} ).

In 1904, Barnes defined multiple gamma functions and multiple zeta
functions, which are given as follows:

In the complex plane, Barnes' multiple zeta function $\zeta
_{r}(s,w|a_{1},\ldots ,a_{r})$ depends on the parameters $a_{1},\ldots
,a_{r} $ (which will be taken nonzero). This function can be defined by the
series%
\begin{equation*}
\zeta _{r}(s,w|a_{1},\ldots ,a_{r})=\sum_{m_{1},\ldots ,m_{r}=0}^{\infty }%
\frac{1}{(w+m_{1}a_{1}+\cdots +m_{r}a_{r})^{s}},
\end{equation*}%
\noindent for $Re(w)>0$ and $Re(s)>r$ \cite{Barnes}.

In \cite{S. N. M. Ruijsenaars}, Ruijsenaars showed how various
known results concerning the Barnes multiple zeta and gamma
functions can be obtained as specializations of simple features
shared by a quite extensive class of functions. The pertinent
functions involve Laplace transforms, and their asymptotic was
obtained by exploiting this. He demonstrated how Barnes' multiple
zeta and gamma functions fit into a recently developed theory of
minimal solutions to first order analytic difference equations.
Both of these approaches to the Barnes functions gave rise to
novel integral representations. In \cite{Kim10}, T Kim studied on
the multiple $L$-series and functional equation of this functions.
He found the value of this function at negative integers in terms
of generalized Bernoulli numbers. In \cite{T. M. Rassia and H. M.
Srivastava}, Russias and Srivastava presented a systematic
investigation of several families of infinite series which are
associated with the Riemann zeta functions, the Digamma
functions,the harmonic numbers, and the stirling numbers of the
first kind. In \cite{K. Matsumoto}, Matsumoto considered general
multiple zeta functions of multi-variables, involving both Barnes
multiple zeta functions and Euler-Zagier sums as special cases. He
proved the meromorphic continuation to the whole space ,
asymptotic expansions, and upper bounded estimates.
These results were expected to have applications to some arithmetical $L$%
-functions. His method was based on \ the classical Mellin-Barnes
integral formula. Ota\cite{K. Ota} studied on Kummer-type
congruences for derivatives of Barnes' multiple Bernoulli
polynomials. Ota generalized these congruences to derivatives of
Barnes' multiple Bernoulli polynomials by an elementary method and
gave a $p$-adic interpolation of them. The Barnes' multiple
Bernoulli polynomials $B_{n}(x,r|a_{1},\ldots ,a_{r})$
are defined by \cite{Barnes}%
\begin{equation}
\frac{t^{r}e^{xt}}{\prod\limits_{j=1}^{r}\left( e^{a_{j}t}-1\right) }%
=\sum_{n=0}^{\infty }B_{n}(x,r|a_{1},\ldots ,a_{r})\frac{t^{n}}{n!},
\label{1.1}
\end{equation}%
\noindent for $|t|<1$. In \cite{Kim9}, T. Kim constructed the Barnes-type
multiple Frobenius-Euler polynomials. He gave a Witt-type formula for these
polynomials. To give the Witt-type formula for Barnes-type multiple
Frobenius-Euler polynomials, he employed the $p$-adic Euler integrals on $%
\mathbb{Z}_{p}$ \cite{Kim5}. He also investigated the properties of the $p$%
-adic Stieltjes transform and $p$-adic Mellin transform. He defined multiple
zeta functions (the Euler-Barnes multiple zeta functions) depending on the
parameters $a_{1},a_{2},\ldots ,a_{r}$ that are taken positive in the
complex number field (\cite{Kim10}, \cite{Kim11}, \cite{Kim12},
\cite{Kim13}%
, \cite{Kim14}, \cite{Kim-Rim}):%
\begin{equation}
\zeta _{r}(s,w,u|a_{1},\ldots ,a_{r})=\sum_{m_{1},\ldots ,m_{r}=0}^{\infty
}%
\frac{u^{-(m_{1}+m_{2}+\ldots +m_{r})}}{(w+m_{1}a_{1}+\cdots
+m_{r}a_{r})^{s}%
},  \label{1.2}
\end{equation}%
\noindent where $Re(w)>0,u\in \mathbb{C}$ with $|u|>1$. These values have a
certain connection with topology and physics, together with the algebraic
relations among them. He also showed that multiple zeta functions can be
continued analytically to $\mathbb{C}^{k}$.

The Euler numbers $E_{n}$ are defined by \cite{H. Tsumura-1},
\begin{equation*}
\frac{2}{e^{t}-1}=\sum_{n=0}^{\infty }E_{n}\frac{t^{n}}{n!}.
\end{equation*}%
They are classical and important in number theory. Frobenius extended
$E_{n}$
to Euler numbers $H^{n}(u)$ belonging to an algebraic number $u$, with $\mid
u\mid >1$, and many authors investigate their properties \cite{Kim4}, \cite%
{Kim9}, . Shiratani and Yamamoto \cite{K. Shiratani and S. Yamamoto}
constructed a $p$-adic interpolation $G_{p}(s,u)$ of the Euler numbers $%
H^{n}(u)$ and as its application, they obtained an explicit formula for $%
L_{p}^{^{\prime }}(0,\chi )$ with any Dirichlet character $\chi $. In \cite%
{H. Tsumura-1}, Tsumura defined the generalized Euler numbers $H_{\chi
}^{n}(u)$ for any Dirichlet character $\chi $, which are analogous to the
generalized Bernoulli numbers. He constructed their $p$-adic interpolation,
which is an extension of Shiratani and Yamamoto's $p$-adic interpolation $%
G_{p}(s,u)$ of $H^{n}(u)$

For $u\in \mathbb{C}$ with $|u|>1$, the Frobenius-Euler polynomial were also
defined by%
\begin{equation}
\frac{1-u}{e^{t}-u}e^{xt}=e^{H(x,u)}=\sum_{n=0}^{\infty }H_{n}(x,u)\frac{%
t^{n}}{n!},  \label{1.3}
\end{equation}%
\noindent where we use the notation by symbolically replacing $H^{m}(x,u)$
by $H_{m}(x,u)$. In the case of $x=0$, the Frobenius-Euler polynomials are
called Frobenius-Euler numbers. We write $H_{m}(u)=H_{m}(0,u)$ (\cite{Kim4},
\cite{Kim9}). Note that $H_{m}(-1)=E_{m}$.

The Frobenius-Euler polynomials of order $r$, denoted by $H_{n}^{(r)}(u,x)$,
were defined as%
\begin{equation}
\left( \frac{1-u}{e^{t}-u}\right) ^{r}e^{xt}=\sum_{n=0}^{\infty
}H_{n}^{(r)}(u,x)\frac{t^{n}}{n!}  \label{1.4}
\end{equation}%
\noindent (\cite{Kim4}, \cite{Kim9}). The values at $x=0$ are called
Frobenius-Euler numbers of order $r$, when $r=1$, the polynomials or numbers
are called ordinary Frobenius-Euler polynomials or numbers. When $x=0$ or $%
r=1$, we often suppress that part of the notation; e.g., $H_{n}^{(r)}(u)$
denotes $H_{n}^{(r)}(u,0)$, $H_{n}(u)$ denotes $H_{n}^{(1)}(u,0)$
\cite{Kim4}%
, \cite{Kim9}. Let $w,a_{1},a_{2},\ldots ,a_{r}$ be complex numbers such
that $a_{i}\neq 0$ for each $i$, $i=1,2,\ldots ,r$. Then the Euler-Barnes'
polynomials of $w$ with parameters $a_{1},\ldots ,a_{r}$ are defined as%
\begin{equation}
\frac{(1-u)^{r}}{\prod\limits_{j=1}^{r}\left( e^{a_{j}t}-u\right) }%
=\sum_{n=0}^{\infty }H_{n}^{(r)}(w,u|a_{1},\ldots ,a_{r})\frac{t^{n}}{n!},
\label{1.5}
\end{equation}%
\noindent for $u\in \mathbb{C}$ with $|u|>1$ \cite{Kim4}, \cite{Kim9}.

In the special case $w=0$, the above polynomials are called the $r$-th Euler
Barnes' numbers. We write%
\begin{equation*}
H_{n}^{(r)}(u|a_{1},\ldots ,a_{r})=H_{n}^{(r)}(0,u|a_{1},\ldots ,a_{r}).
\end{equation*}

In \cite{Kim4}, \cite{Kim9}, by using $p$-adic (Euler) integrals on
$\mathbb{%
Z}_{p}$, Kim and Rim constructed Changhee-Barnes' $q$-Euler numbers and
polynomials which are related to the $q$-analogue of Euler-Barnes'
polynomials and numbers. The $q$-analogue of Frobenius-Euler numbers, by
using $p$-adic Euler integral, are given as follows \cite{Kim4},
\cite{Kim9}:%
\begin{equation*}
H_{n,q}(u)=\int\limits_{\mathbb{Z}_{p}}[x]_{q}^{n}d\mu _{u}(x)=\frac{1-u}{%
(1-q)^{n}}\sum_{l=0}^{n}\binom{n}{l}(-1)^{l}\frac{1}{1-uq^{l}},
\end{equation*}%
\noindent where $p$-adic $q$-integral is given by: Let $f$ be uniformly
differentiable function at a point of $\mathbb{Z}_{p}$.
Then we have \cite{Kim-Rim}%
\begin{equation*}
\int\limits_{\mathbb{Z}_{p}}f(x)d\mu _{u}(x)=\underset{N\rightarrow \infty
}{%
lim}[p^{N}:u]^{-1}\sum_{0\leq x<p^{N}}f(x)u^{x},
\end{equation*}%
\noindent with $u\in \mathbb{C}_{p}$, where $\mathbb{C}_{p}$ is the
completion of algebraic closure of $\mathbb{Q}_{p}$. In the above we assume
that $u\in \mathbb{C}_{p}$ with $|1-u|_{p}\geq 1$. If $q\in \mathbb{C}_{p}$,
one normally assumes that $|q|<1$ and $|1-q|_{p}\leq p^{-1/(p-1)}$, so that
$%
q^{x}=exp(x\log q)$. We use the notation%
\begin{equation*}
\lbrack x]=[x]_{q}=\frac{1-q^{x}}{1-q}\text{, and
}[x:z]_{q}=\frac{1-z^{x}}{%
1-z}
\end{equation*}%
\noindent (\cite{Kim7}, \cite{Kim9}, \cite{Kim-Rim}). The $q$-analogue of
Euler-Barnes' multiple numbers, which reduce to Euler-Barnes' multiple
numbers at $q=1$, are given as follows \cite{Kim-Rim}:

Let $a_{1},\ldots ,a_{r}$; $b_{1},\ldots ,b_{r}$ will be taken in the
nonzero integers and let $q\in \mathbb{C}_{p}$ with $|1-q|_{p}<1,w\in
\mathbb{Z}_{p},u\in \mathbb{C}_{p}$ with $|1-u|_{p}\geq 1$. Then \cite%
{Kim-Rim}%
\begin{eqnarray}
&&H_{n,q}^{(r)}(u,w|a_{1},\ldots ,a_{r}:b_{1},\ldots ,b_{r})  \notag \\
&=&\int\limits_{\mathbb{Z}_{p}^{r}}[w+\sum_{j=1}^{r}a_{j}x^{j}]^{n}u^{%
\sum_{j=1}^{r}(b_{j}-1)x_{j}}d\mu _{u}(x_{1})\cdots \mu _{u}(x_{r})
\label{1.6}
\end{eqnarray}

\begin{theorem}
\label{th1} For $n\geq 0$, we have%
\begin{eqnarray}
&&H_{n,q}^{(r)}(u,w|a_{1},\ldots ,a_{r}:b_{1},\ldots ,b_{r})  \notag \\
&=&\frac{(1-u)^{r}}{(1-q)^{n}}\sum_{l=0}^{n}\binom{n}{l}(-1)^{l}q^{lw}\frac{1%
}{\prod\limits_{j=1}^{r}\left( 1-q^{la_{j}}u^{b_{j}}\right) },
\text{ cf. \cite{Kim-Rim}}. \label{1.7}
\end{eqnarray}
\end{theorem}

\noindent Note that%
\begin{equation*}
\underset{q\rightarrow 1}{\text{lim}}H_{n,q}^{(r)}(u,w|a_{1},\ldots
,a_{r}:1,\ldots ,1)=H_{n}^{(r)}(u^{-1},w|a_{1},\ldots ,a_{r}).
\end{equation*}

In this paper, we define generating functions of Changhee-Barnes'
$q$-Euler numbers and polynomials on $\mathbb{C}$-plane. These
functions are very important to construct multiple zeta functions.
By using Mellin transform's formula and Cauchy Theorem , we prove
the analytic continuation of Barnes' type multiple Changhee
$q$-zeta function. Finally we give some relations between Barnes'
type multiple Changhee $q$-zeta function and Barnes' type
multiple\ generalized Changhee $q$-Bernoulli numbers.
 Moreover, we give the value of these functions at negative
integers.

\section{Changhee-Barnes' Type $q$-Euler Numbers and Polynomials}

In this chapter, we define generating function of the Changhee-Barnes' type
$%
q$-Euler numbers and polynomials as follows:

Let $u$ be the algebraic element of the complex number field $\mathbb{C}$,
with $|u|<1$. For $w_{1},v_{1},s\in \mathbb{C}$ with $%
Re(w_{1})>0,Re(v_{1})>0 $, we then define%
\begin{eqnarray}
F_{u^{-1},q}(t|w_{1};v_{1}) &=&(1-u)e^{t/(1-q)}\sum_{j=0}^{\infty
}\frac{1}{%
(1-q)^{j}}(-1)^{j}\frac{1}{1-q^{w_{1}j}u^{v_{1}}}\frac{t^{j}}{j!}  \notag \\
&=&\sum_{k=0}^{\infty }H_{k,q}(u^{-1}|w_{1};v_{1})\frac{t^{k}}{k!},\text{ }%
|t|<2\pi .  \label{2.1}
\end{eqnarray}

\begin{theorem}
\label{th2} Let $u$ be algebraic integer with $|u|<1,u\in \mathbb{C}$ and $%
w_{1},v_{1}\in \mathbb{C}$ with $Re(w_{1})>0,Re(v_{1})>0$. Then we have%
\begin{equation*}
H_{k,q}(u^{-1}|w_{1};v_{1})=\frac{1-u}{(1-q)^{k}}\sum_{j=0}^{\infty
}\binom{k%
}{j}(-1)^{j}\frac{1}{1-q^{w_{1}j}u^{v_{1}}}.
\end{equation*}
\end{theorem}

\begin{proof}
By using (\ref{2.1}) we have%
\begin{eqnarray*}
&&\sum_{k=0}^{\infty }H_{k,q}(u^{-1}|w_{1};v_{1})\frac{t^{k}}{k!} \\
\qquad \qquad &=&(1-u)\left( \sum_{n=0}^{\infty }\frac{1}{(1-q)^{n}}\frac{%
t^{n}}{n!}\right) \left( \sum_{j=0}^{\infty }\frac{1}{(1-q)^{j}}(-1)^{j}%
\frac{1}{1-q^{w_{1}j}u^{v_{1}}}\frac{t^{j}}{j!}\right) .
\end{eqnarray*}%
\noindent By applying Cauchy product in the above, we easily obtain%
\begin{equation*}
\sum_{k=0}^{\infty }H_{k,q}(u^{-1}|w_{1};v_{1})\frac{t^{k}}{k!}%
=(1-u)\sum_{k=0}^{\infty }\left(
\sum_{j=0}^{k}\binom{k}{j}(-1)^{j}\frac{1}{%
1-q^{w_{1}j}u^{v_{1}}}\right) \frac{1}{(1-q)^{k}}\frac{t^{k}}{k!}.
\end{equation*}%
\noindent Now, by comparing coefficients of $\frac{t^{k}}{k!}$ in the above
series, we arrive at the desired result.
\end{proof}

\begin{remark}
\label{re1} We note that if $v_{1}=1$ in (\ref{2.1}), then (\ref{2.1})
reduces to (8.1) in \cite{Kim-Simsek}. By substituting $r=1$ and $%
v_{1}=v_{2}=\cdots =v_{r}=1$ in (\ref{1.6}) and (\ref{1.7}), then Changhee
$%
q $-Euler numbers reduce to Daehee numbers. Indeed, by using (\ref{2.1}), we
obtain%
\begin{eqnarray*}
F_{u^{-1},q}(t|w_{1};v_{1}) &=&(1-u)\sum_{l=0}^{\infty
}u^{v_{1}l}e^{t/(1-q)}\sum_{j=0}^{\infty }\frac{1}{(1-q)^{j}}%
(-1)^{j}q^{w_{1}jl}\frac{t^{j}}{j!} \\
&=&(1-u)\sum_{l=0}^{\infty }u^{v_{1}l}e^{t[w_{1}l]_{q}}.
\end{eqnarray*}%
\noindent Hence%
\begin{eqnarray}
F_{u^{-1},q}(t|w_{1};v_{1}) &=&(1-u)\sum_{l=0}^{\infty
}u^{v_{1}l}e^{t[w_{1}l]_{q}}  \notag \\
&=&\sum_{n=0}^{\infty }H_{n,q}(u^{-1}|w_{1};v_{1})\frac{t^{n}}{n!}.
\label{2.2}
\end{eqnarray}%
\noindent In the above for $v_{1}=1$ the function $F_{u^{-1},q}(t|w_{1};1)$
is the generating function of $q$-Daehee polynomials.
\end{remark}

We now define the generating function of Changhee $q$-Euler polynomials as
follows:%
\begin{eqnarray}
F_{u^{-1},q}(t,w|w_{1};v_{1})
&=&e^{[w]_{q}t}F_{u^{-1},q}(q^{w}t|w_{1};v_{1})
\notag \\
&=&\sum_{n=0}^{\infty }H_{n,q}(u^{-1},w|w_{1};v_{1})\frac{t^{n}}{n!},\text{
}%
|t|<2\pi .  \label{2.3}
\end{eqnarray}%
\noindent By using (\ref{2.2}) and (\ref{2.3}), we easily see that%
\begin{eqnarray}
F_{u^{-1},q}(t,w|w_{1};v_{1}) &=&e^{[w]_{q}t}(1-u)\sum_{l=0}^{\infty
}u^{v_{1}l}e^{[w_{1}l]_{q}q^{w}t}  \notag \\
&=&(1-u)\sum_{l=0}^{\infty }u^{v_{1}l}e^{[w+w_{1}l]_{q}t}.  \label{2.4}
\end{eqnarray}

\noindent By (\ref{2.3}) and (\ref{2.4}), we obtain%
\begin{eqnarray*}
F_{u^{-1},q}(t,w|w_{1};v_{1}) &=&(1-u)\sum_{l=0}^{\infty
}u^{v_{1}l}e^{[w+w_{1}l]_{q}t} \\
&=&\sum_{n=0}^{\infty }H_{n,q}(u^{-1},w|w_{1};v_{1})\frac{t^{n}}{n!}.
\end{eqnarray*}

Now, we will generalize (\ref{2.1}). Let $w_{1},w_{2},\ldots
,w_{r};v_{1},v_{2},\ldots ,v_{r}\in \mathbb{C}$. Then generating function of
multiple Changhee $q$-Euler polynomials are given as follows:%
\begin{eqnarray}
&&F_{u^{-1},q}(t,w|w_{1},w_{2},\ldots ,w_{r};v_{1},v_{2},\ldots ,v_{r})
\notag \\
&=&(1-u)\sum_{n_{1},n_{2},\ldots ,n_{r}=0}^{\infty
}u^{\sum_{i=1}^{r}n_{i}v_{i}}e^{\left[ \sum_{i=1}^{r}n_{i}w_{i}\right]
_{q}t}
\notag \\
&=&\sum_{n=0}^{\infty }H_{n,q}^{(r)}(u^{-1},w|w_{1},w_{2},\ldots
,w_{r};v_{1},v_{2},\ldots ,v_{r})\frac{t^{n}}{n!},\text{ }|t|<2\pi .
\label{2.5}
\end{eqnarray}

\noindent By using (\ref{2.5}), we arrive at the following theorem easily:

\begin{theorem}
\label{th3} For $n\geq 0$, we have%
\begin{eqnarray*}
&&H_{n,q}^{(r)}(u^{-1},w|w_{1},w_{2},\ldots ,w_{r};v_{1},v_{2},\ldots
,v_{r})
\\
&=&\frac{(1-u)^{r}}{(1-q)^{n}}\sum_{l=0}^{n}\binom{n}{l}(-1)^{l}q^{lw}\frac{1%
}{\prod\limits_{j=1}^{r}\left( 1-q^{la_{j}}u^{b_{j}}\right) }.
\end{eqnarray*}
\end{theorem}

\noindent The proof of this theorem is also given by Kim and Rim \cite%
{Kim-Rim}. Their proof is related to $p$-adic Euler integration.

Note that%
\begin{equation*}
\underset{q\rightarrow 1}{\text{lim}}H_{n,q}^{(r)}(u^{-1},w|w_{1},w_{2},%
\ldots ,w_{r}:1,\ldots ,1)=H_{n}^{(r)}(u^{-1},w|w_{1},w_{2},\ldots ,w_{r}).
\end{equation*}

\section{Analytic Continuation of Changhee $q$-Euler-Zeta Function}

By using (\ref{2.1}) and (\ref{2.5}), we consider Changhee
$q$-zeta functions \noindent $\zeta _{q}^{(r)}(u,w|w_{1},\ldots
,w_{r};v_{1},\ldots ,v_{r})$. We also give analytic continuation
of \noindent $\zeta _{q}^{(r)}(u,w|w_{1},\ldots
,w_{r};v_{1},\ldots ,v_{r})$ by using Mellin transform's formula
of (\ref{2.1}) and (\ref{2.5}). Let $w,w_{1,}v_{1}$ be complex
numbers with positive real parts. For $s\in
\mathbb{C}$, we consider Changhee $q$-zeta functions as follows:%
\begin{equation}
\frac{1}{1-u}\frac{1}{\Gamma (s)}\int\limits_{0}^{\infty
}F_{u^{-1},q}(-t,w|w_{1};v_{1})t^{s-1}dt=\sum_{n=1}^{\infty
}\frac{u^{v_{1}n}%
}{[w+w_{1}n]_{q}^{s}}.  \label{3.1}
\end{equation}

\noindent ({\ref{3.1}) is the Mellin transform of (\ref{2.1}). By (%
\ref{3.1}), we define Changhee $q$-zeta function as follows:

For $s\in \mathbb{C}$,%
\begin{equation*}
\zeta _{q}(s,w,u|w_{1};v_{1})=\sum_{n=0}^{\infty }\frac{u^{v_{1}n}}{%
[w+w_{1}n]_{q}^{s}}.
\end{equation*}

\noindent Note that, for $v_{1}=1$, $\zeta _{q}(s,w,u|w_{1};1)$ reduces to %
\cite{Kim-Simsek}. Thus, $\zeta _{q}(s,w,u|w_{1};v_{1})$ is the analytic
continuation on $\mathbb{C}$, with simple pole at $s=1$.

By using Cauchy Theorem and Residue Theorem in (\ref{3.1}), for positive
integer $n$, we easily arrive at the following theorem:

\begin{theorem}
\label{th4} For $n\in \mathbb{Z}^{+}$,%
\begin{equation*}
\zeta _{q}(-n,w,u|w_{1};v_{1})=\frac{1}{1-u}H_{n,q}(u^{-1},w|w_{1};v_{1}).
\end{equation*}
\end{theorem}

Similarly we define multiple Changhee $q$-zeta function as
follows:

Let $w,w_{1},w_{2},\ldots ,w_{r},v_{1},v_{2},\ldots ,v_{r}$ be
complex numbers with positive real parts and $r\in
\mathbb{Z}^{+}$. Then we construct multiple Changhee $q$-zeta
functions as follows: For $s\in
\mathbb{C}$,%
\begin{eqnarray}
&&\frac{1}{(1-u)^{r}}\frac{1}{\Gamma (s)}\int\limits_{0}^{\infty
}F_{u^{-1},q}^{(r)}(t,w|w_{1},w_{2},\ldots ,w_{r};v_{1},v_{2},\ldots
,v_{r})t^{s-1}dt  \notag \\
&=&\sum_{n_{1},n_{2},\ldots ,n_{r}=0}^{\infty }\frac{u^{%
\sum_{i=1}^{r}v_{i}n_{i}}}{\left[ w+\sum_{i=1}^{r}w_{i}n_{i}\right]
_{q}^{s}}%
.  \label{3.2.2}
\end{eqnarray}

\noindent The above equation is the Mellin transform of
(\ref{2.5}). By using (\ref{3.2.2}), we define multiple Changhee
$q$-zeta functions as
follows:%
\begin{equation*}
\zeta _{q}^{(r)}(s,u,w|w_{1},\ldots ,w_{r};v_{1},\ldots
,v_{r})=\sum_{n_{1},n_{2},\ldots ,n_{r}=0}^{\infty }\frac{%
u^{\sum_{i=1}^{r}v_{i}n_{i}}}{\left[ w+n_{1}w_{1}+\cdots +n_{r}w_{r}\right]
_{q}^{s}}.
\end{equation*}

\noindent Note that, for $v_{1}=v_{2}=\cdots =v_{r}=1$, $\zeta
_{q}^{(r)}(s,u,w|w_{1},\ldots ,w_{r};1,\ldots ,1)$ is reduced to
Definition 6 in \cite{Kim-Simsek}. By using Cauchy Theorem and
Residue Theorem in (\ref{3.2.2}), we easily obtain the following
theorem:

\begin{theorem}
\label{th5} For $n\in \mathbb{Z}^{+}$, we have%
\begin{equation*}
\zeta _{q}^{(r)}(-n,w,u|w_{1},\ldots ,w_{r};v_{1},\ldots ,v_{r})=\frac{1}{%
(1-u)^{r}}H_{n,q}(u^{-1},w|w_{1},\ldots ,w_{r};v_{1},\ldots ,v_{r}).
\end{equation*}
\end{theorem}


\begin{thebibliography}{99}
\bibitem{Andrews} G. E. Andrews, $q$-analogs of the binomal coefficient
congruences of Babbage, Wolstenhome and Glaiser, Discrete Math., 204 (1999),
15-25.

\bibitem{Askey} R. Askey, The $q$-gamma and $q$-beta functions, Appl. Anal.,
8 (1978), 125-141.

\bibitem{Barnes} W. Barnes, On theory of the multiple gamma functions,
Trans. Camb. Philos. Soc., 19 (1904), 374-425.

\bibitem{Cherednik} I. Cherednik, On $q$-analogues of the Riemann's zeta
function, Selecta Math., Vol. 7\ (5) (2001), 447-491.

\bibitem{E. Friedman and S. Ruijsenaars} E. Friedman and S. Ruijsenaars,
Shintani-Barnes zeta and gamma functions, Advences in Math., 187 (2004),
362-395.

\bibitem{M. Jimbo and T. Miwa} M. Jimbo and T. Miwa, Quantum KZ equation
with $\mid q\mid =1$ and correlation functions of the $XXZ$ model in the
gapless regime, J. Phys. A, 29 (1996), 2923-2958.

\bibitem{M. Nishizawa} M. Nishizawa, On a $q$-analogue of the multiple gamma
functions, Lett. Math. Phys., 37 (1996), 201-209.

\bibitem{Katriel} J. Katriel, Stirling numbers identities interconsistency
of $q$-analogues, J. Phys. A: Gen. Math. 31 (1988), 3559-3572.

\bibitem{Kim1} T. Kim, On explicit formulas of $p$-adic $q$-$L$-functions,
Kyushu J. Math., Vol. 48 (1) (1994), 73-86.

\bibitem{Kim2} T. Kim, On a $q$-analogue of the $p$-adic Log gamma functions
and related integrals, J. Number Theo., 76 (1999), 320-329.

\bibitem{Kim4} T. Kim, An invariant $p$-adic integral associated with Daehee
Numbers, Integral Transform and Special Functions, 13 (2002), 65-69.

\bibitem{Kim5} T. Kim, $q$-Volkenborn Integration, Russ. J. Math. Phys., 19
(2002), 288-299.

\bibitem{Kim7} T. Kim, Non-archimedean $q$-integrals associated with
multiple Changhee $q$-Bernoulli Polynomials, Russ. J. Math Phys., Vol. 10
(2003), 91-98.

\bibitem{Kim9} T. Kim, On Euller-Barnes multiple zeta functions, Russ. J.
Math Phys., Vol. 10 (3) (2003), 261-267.

\bibitem{Kim10} T. Kim, A note on Dirichlet series, Proc. Jangjeon Math.
Soc., 6(2) (2003),161-166.

\bibitem{Kim11} T. Kim, Some of powers of consecutive $q$-integrals, Adv.
Stud. Contep. Math., 9 (1) (2004), 15-18.

\bibitem{Kim12} T. Kim, A note on $q$-zeta functions, Proceedings of The
15th International Conference of The Jangjeon Mathematical Society, August
5-7, S. Korea (2004), 110-114.

\bibitem{Kim13} T. Kim, A note on the $q$-multiple zeta function, Adv. Stud.
Contep. Math., 8 (12) (2004), 111-113.

\bibitem{Kim14} T. Kim, $p$-adic $q$-integrals associated with the
Changhee-Barnes' $q$-Bernoulli Polynomials, Integral Transform and Special
Functions, Vol, Month (2004), 1-6.

\bibitem{Kim-Rim} T. Kim and S. -H. Rim, A note on the $q$-integrals and
$q$%
-series, Adv. Stud. Contep. Math., 2 (2000), 37-45.

\bibitem{T.Kim-S.-H.Rim} T. Kim and S. -H. Rim, On Changhee-Barnes' $q$%
-Euler numbers and polynomials, Adv. Stud. Contep. Math., 9 (2) (2004),
81-86.

\bibitem{Kim-Simsek} T. Kim, Y. Simsek, H. M. Srivastava, $q$-Bernoulli
Numbers and Polynomials Associated with Multiple $q$-Zeta Functions and
Basic $L$-series, Submitted.

\bibitem{Khrennikov} A. Khrennikov, $p$-Adic Valued Distributions in
Mathematical Physics, Kluwer Academic Publishers, London (1994).

\bibitem{T. H. Koornwinder} T. H. Koornwinder, Special functions and $q$%
-commuting valuables, Fields Institute Communications Vol. 14 (1997).

\bibitem{N. Kurokawa} N. Kurokawa, Multiple sine functions and selberg zeta
functions, Proc. Japan Acad. A, 67 (1991), 61-64.

\bibitem{K. Matsumoto} K. Matsumoto, The analytic continuation and the
asyptotic behaviour of certain multiple zeta-function I, J. Number Theory
101 (2003), 223-243.

\bibitem{K. Ota} K. Ota, On Kummer-type congruences for derivatives of
Barnes' multiple Bernoulli Polynomials, J. Number Theory, 92 (2002), 1-36.

\bibitem{C. A. Nelson and M. G. Gartley} C. A. Nelson and M. G. Gartley, On
the zeros of the $q$-analogue exponential function, J. Phys. A: Gen. Math.,
27 (1994), 3857-3881.

\bibitem{T. M. Rassia and H. M. Srivastava} T. M. Rassia and H. M.
Srivastava, Some classes of infinite series associated with the Riemann zeta
and polygamma functions and generalized harmonic numbers, Appl. Math.
Computation, 131 (2002), 593-605.

\bibitem{A. M. Robert} A. M. Robert, A course in $p$-adic Analysis, Springer
(2000).

\bibitem{S. N. M. Ruijsenaars} S. N. M. Ruijsenaars, On Barnes' multiple
zeta function and gamma functions, Adv. in Math., 156, (2000), 107-132.

\bibitem{B. E. Sagan} B. E. Sagan, Congruence properties of $q$-analogs,
Adv. in Math., 95, (1992), 127-143.

\bibitem{W. H. Schikhof} W. H. Schikhof, Ultrametric Calculus, Cambridge
Univ. Press (1984).

\bibitem{K. Shiratani} K. Shiratani, On Euler numbers, Mem. Fac. Kyushu
Uni., 27 (1) (1973), 1-5.

\bibitem{K. Shiratani and S. Yamamoto} K. Shiratani and S. Yamamoto, On a
$p$%
-adic interpolation function for the Euler numbers and its derivative, Mem.
Fac. Kyushu Uni., 39 (1985), 113-125.

\bibitem{H. Tsumura-1} H. Tsumura, On a $p$-adic interpolation of
generalized Euler Numbers and its applications, Tokyo J. Math., 10 (2)
(1987), 281-293.

\bibitem{Washington} Washington, Introduction to Cyclomotic field, GTM 83,
Springer-Verlag (1982).

\bibitem{E. T. Wittaker and G. N. Watson} E. T. Wittaker and G. N. Watson, A
course of modern Analysis, Cambridge Univ. Press (1990).
\end{thebibliography}
\end{document}